\documentclass{gen-j-l}

\usepackage{amsmath,amssymb,amsfonts,epic,longtable,curves}

\newtheorem{theorem}{Theorem}[section]
\newtheorem{lemma}[theorem]{Lemma}
\newtheorem{proposition}[theorem]{Proposition}

\theoremstyle{definition}

\theoremstyle{remark}

\numberwithin{equation}{section}

\def\DJ{{\hbox{D\kern-.8em\raise.15ex\hbox{--}\kern.35em}}}

\def\NSERC{The first author was supported in part by the
NSERC Grant A-5285.}

\def\bC{{\mbox{\bf C}}}
\def\sA{{\mathcal A}}
\def\GL{{\mbox{\rm GL}}}
\def\SL{{\mbox{\rm SL}}}
\def\ve{{\varepsilon}}
\def\diag{{\mbox{\rm diag}}}

\begin{document}

\title[Hermitian and skew-hermitian matrices]{Solution of the
congruence problem for arbitrary hermitian and
skew-hermitian matrices over polynomial rings}

\author[D.\v Z. \DJ okovi\'c]{Dragomir \v Z. \DJ okovi\'c}
\address{Department of Pure Mathematics, University of Waterloo,
Waterloo, Ontario, N2L 3G1, Canada}
\thanks{\NSERC\newline e-mail: djokovic@uwaterloo.ca,
fszechtm@herod.uwaterloo.ca}

\author[F. Szechtman]{Fernando Szechtman}

\date{}

\begin{abstract}
Let ${}^*$ be the involutorial automorphism of the
complex polynomial algebra $\bC[t]$ which sends $t$ to $-t$.
Answering a question raised by V.G. Kac,
we show that every hermitian or skew-hermitian matrix
over this algebra is congruent to the direct sum of
$1\times1$ matrices and $2\times2$ matrices with zero diagonal.
Moreover we show that if two $n\times n$ hermitian or
skew-hermitian matrices have the same
invariant factors, then they are congruent.
The complex field can be replaced by any algebraically closed field
of characteristic $\ne2$.
\end{abstract}

\maketitle

\section{Introduction}

Let $R$ be the polynomial algebra $F[t]$ in one variable $t$
over an algebraically closed field $F$ of characteristic $\ne2$.
Let ${}^*$ denote the involution of $R$
which is the identity on $F$ and sends $t$ to $-t$.
(We remark that all nontrivial $F$-involutions of $F[t]$
are conjugate in ${\rm Aut}_F(F[t])$.)
It induces the $Z_2$-gradation $R=R_0\oplus R_1$ of $R$ with
$R_0=F[t^2]$ and $R_1=tR_0$. We shall refer to the
elements of $R_0$ (resp. $R_1$) as {\em even} (resp. {\em odd}).

Let $M_n(R)$ denote the algebra of $n$ by $n$ matrices over $R$.
If $A=(a_{ij})\in M_n(R)$, we define $A^*$ to be the matrix
$B=(b_{ij})\in M_n(R)$ where $b_{ij}=a_{ji}^*$. Thus ${}^*$ is
now made into an involution of $M_n(R)$.
We say that $A\in M_n(R)$ is {\em hermitian}
(resp. {\em skew-hermitian}) if $A^*=A$ (resp. $A^*=-A$).
Two hermitian (resp. skew-hermitian) matrices $A,B\in M_n(R)$
are said to be {\em congruent} if $B=S^*AS$ for some
$S\in\GL_n(R)$.

Not long ago V.G. Kac \cite{VK} posed the following question to
the first author (see also \cite{BKL}).

{\em If $F$ is the complex field, is it true that every hermitian
or skew-hermitian matrix $A\in M_n(R)$ is congruent to the
direct sum of $1\times1$ matrices and $2\times2$ matrices
with zero diagonal?}

Note that no condition is imposed on the determinant of $A$.
(The usual restriction is that $A$ be unimodular.)
The two cases, hermitian and skew-hermitian, of this
problem are tightly linked because if $A$ is hermitian
then $tA$ is skew-hermitian, and vice versa.

The main objective of our paper is to give an affirmative
answer to Kac's question (Theorem \ref{glavni}),
which we find quite surprising.
The case $n=2$ is dealt with in Section 3 and
the general case in Section 4.
In Section 4, we also prove that two hermitian (or skew-hermitian)
matrices $A,B\in M_n(R)$ are congruent if and only if 
they have the same invariant factors (Theorem \ref{invar}).
Then, in Section 5, we are able to characterize the sequence
of invariant factors of a hermitian or skew-hermitian matrix,
and to give the canonical form under congruence for such a matrix.
In the last section we make comments on other fields and
characterize those for which Kac's question has positive answer.

The authors would like to thank Prof. Kac for his interest in
our work and for proposing this interesting problem.
We also thank Prof. L. Vaserstein for his comments on an
earlier version of this paper.

\section{Preliminaries}

The elements $a\in R$ are polynomials and so they can be
evaluated at any point $\lambda\in F$. We denote
by $a(\lambda)$ the value of $a$ at $\lambda$. 
We say that a nonzero element $a\in R$ is {\em pure} if
$\gcd(a,a^*)=1$.
If $a,b\in R$ with $a$ pure and $b$ even (resp. odd),
then there exists $x\in R$ such that $ax+a^*x^*=b$
(resp. $ax-a^*x^*=b$). (Choose $y,z\in R$ such that
$ay+a^*z=b/2$ and set $x=y+z^*$ (resp. $x=y-z^*$).)

If $I=Ra$ is a homogeneous (i.e., $*$-invariant) ideal of $R$,
then its generator $a$ is also homogeneous, i.e., it is either
even or odd. If $A=(a_{ij})\in M_n(R)$ is hermitian or
skew-hermitian, then the ideal generated by all entries
$a_{ij}$ is $*$-invariant and we denote its generator
by $\gcd(A)$. Hence $\gcd(A)$ is the first invariant
factor of $A$.

Let us fix a hermitian or skew-hermitian matrix $A\in M_n(R)$.
Let $R^n$ denote the free $R$-module of rank $n$
consisting of column vectors.
We shall denote by $e_1,\ldots,e_n$ the standard basis vectors
of $R^n$. The matrix $A$ defines a hermitian
or skew-hermitian form $f_A:R^n\times R^n\to R$ by
\begin{equation*}
f_A(v,w)=v^*Aw.
\end{equation*}

By \cite[Lemma 1]{DZ}, $A$ is congruent to the direct sum 
of a zero matrix and a hermitian or skew-hermitian
matrix with nonzero determinant.
(The proof given there in the hermitian case is also
valid in the skew-hermitian case.)
This argument shows that it suffices to
consider only the hermitian or skew-hermitian matrices
with nonzero determinant.

As $F$ is algebraically closed, if $n\ge2$ there exist nonzero
isotropic vectors, i.e., nonzero vectors $v\in R^n$ such that
$f_A(v,v)=0$.

Assume that $\det(A)\ne0$ and set $d=\gcd(A)$.
Then $A=dB$ for some matrix $B\in M_n(R)$ such that
$B^*=\pm B$ and $\gcd(B)=1$.
Therefore, without any loss of generality, we may assume
that $\det(A)\ne0$ and $\gcd(A)=1$.

\section{The case $n=2$}

In this section we show that the answer to Kac's
question is affirmative if $n=2$. We start with the
hermitian case.

\begin{proposition} \label{her-2}
If $A=A^*\in M_2(R)$, $\det(A)\ne0$, and $\gcd(A)=1$,
then $A$ is congruent to $\diag(1,\det(A))$.
\end{proposition}

\begin{proof}
Since there exist nonzero isotropic vectors, we may assume that 
$$
A=\left(\begin{array}{cc}0&a\\a^*&b\end{array}\right).
$$
The element $a_0=\gcd(a,a^*)$ is homogeneous, i.e.,
$a_0^*=\pm a_0$. We have a factorization $a=a_0a_1$ where
$a_1$ is pure.
By the hypothesis, $\gcd(a_0,b)=1$.
Consequently, there exist homogeneous elements $x$ and $y$,
with $y$ even, such that $a_0x+by=1$. Clearly we may
assume that $y(0)\ne0$. Choose a factorization
$y=zz^*$ such that $a_1z$ is pure. Then there exists
$w\in R$ such that
$$
a_1zw+a_1^*z^*w^*=1.
$$
Since $a_0x$ is even, we find that
\begin{eqnarray*}
1 &=& a_0x+by \\
&=& a_0x(a_1zw+a_1^*z^*w^*)+bzz^* \\
&=& axwz+a^*x^*w^*z^*+bzz^* \\
&=& f_A(x^*w^*e_1+ze_2,x^*w^*e_1+ze_2).
\end{eqnarray*}
The assertion of the proposition is now obvious.
\end{proof}

Next we consider the skew-hermitian case. We shall need
the following simple lemma.

\begin{lemma} \label{nzd-lema}
Let $a,b\in R$ satisfy $\gcd(a,b)=\gcd(b,b^*)=1$. Then there exist
$x,y\in R$, with $x$ even, such that $ax+by=1$.
\end{lemma}

\begin{proof}
Choose $u,v\in R$ such that $au+bv=1$  and $z\in R$
such that $bz-b^*z^*=u^*-u$.
Then $x=u+bz\in R_0$ and $y=v-az$ satisfy $ax+by=1$.
\end{proof}

\begin{proposition} \label{sk-2}
Let $A\in M_2(R)$, $A^*=-A$, $\det(A)\ne0$, and $\gcd(A)=1$. Then
$A$ is congruent to a matrix with zero diagonal.
\end{proposition}

\begin{proof}
We may assume that
$$
A=\left(\begin{array}{cc}0&a\\-a^*&b\end{array}\right).
$$
As $a(0)\ne0$, we can write $a=a_1cc^*$ with $a_1c$ pure.
By Lemma \ref{nzd-lema}, there exist $x,d\in R$, with $x$ even,
such that
$$
bx+cd=1.
$$
By replacing $(x,d)$ with $(x+\lambda cc^*,d-\lambda bc^*)$,
where $\lambda\in F$ is suitably chosen, we may
assume that $\gcd(a_1,x)=1$. Since $b$ is odd, we can choose $w\in R$
such that $c^*w^*-cw=b$. Since $\gcd(a_1^*,cx)=1$, there exist
$v,p\in R$ such that
$$
a_1^*v-cxp^*=xw^*-d^*.
$$
Choose $q\in R$ such that
$$
a_1^*q-a_1q^*=p-p^*
$$
and set
$$
y=v+cxq,\quad z=w+c^*(p+a_1q^*).
$$
Then
\begin{eqnarray*}
c^*z^*-cz &=& c^*w^*+cc^*(p^*+a_1^*q)-cw-cc^*(p+a_1q^*) \\
&=& c^*w^*-cw=b, \\
xz^*-a_1^*y &=& xw^*+cx(p^*+a_1^*q)-a_1^*v-a_1^*cxq \\
&=& cxp^*+xw^*-a_1^*v = d^*.
\end{eqnarray*}
Hence
\begin{eqnarray*}
a_1^*c^*y &=& c^*xz^*-c^*d^*=c^*xz^*-(1+bx) \\
&=& x(c^*z^*-b)-1=cxz-1,
\end{eqnarray*}
and so
$$
S=\left(\begin{array}{ll}cx&y\\a_1^*c^*&z\end{array}\right)
\in\SL_2(R).
$$
We have
$$
S^*\left(\begin{array}{cc}0&a_1c^2\\-a_1^*(c^*)^2&0\end{array}
\right)S=\left(\begin{array}{cc}0&r\\-r^*&s\end{array}\right),
$$
where
\begin{eqnarray*}
r &=& a_1cc^*(cxz-a_1^*c^*y)=a_1cc^*=a, \\
s &=& a_1c^2y^*z-a_1^*(c^*)^2yz^* \\
&=& a_1cy^*\cdot cz-a_1^*c^*y\cdot c^*z^* \\
&=& (c^*xz^*-1)cz-(cxz-1)c^*z^* \\
&=& c^*z^*-cz=b.
\end{eqnarray*}
Hence
$$
\left(\begin{array}{cc}0&r\\-r^*&s\end{array}\right)=A.
$$
\end{proof}

\section{Equivalence implies congruence}

The first main theorem is a simple consequence of the following two
propositions. The first one deals with hermitian matrices.

\begin{proposition} \label{her-n}
If $A=(a_{ij})\in M_n(R)$, $A^*=A$ and $\gcd(A)=1$,
then there exists $v\in R^n$ such that $f_A(v,v)=1$.
\end{proposition}

\begin{proof} We may assume that $\det(A)\ne0$.
The proof will be by induction on $n$. The case $n=1$ is obvious.
For the case $n=2$ see Proposition \ref{her-2}.
Thus let $n>2$.

Since there exist nonzero isotropic vectors, we may assume that
$a_{11}=0$. Moreover, we may assume that $a_{1j}=0$ for $j<n$.
Denote by $E_{ij}$ the matrix of order $n$ whose $(i,j)$-th
entry is 1 and all other entries are 0, and by $I_n$ the
identity matrix.
For any $\lambda\in F$, the matrix
$$
A_\lambda=(I_n+\lambda E_{21})A(I_n+\lambda E_{12})
$$
is congruent to $A$.
Let $A'_\lambda$ denote the submatrix of
$A_\lambda$ obtained by deleting the first row and column.
Set $A'=A'_0$. Note that for $\lambda,\mu\in F$ we have
$$
A'_\lambda-A'_\mu=
(\lambda-\mu)(a_{1n}E_{1,n-1}+a_{n1}E_{n-1,1}),
$$
where now $E_{ij}$'s have order $n-1$.

As $\det(A)\ne0$, we have $a_{rs}\ne0$ for some
$r,s\in\{2,3,\ldots,n-1\}$.
Since $a_{rs}$ has only finitely many monic divisors, there
exist $\lambda,\mu\in F$, with $\lambda\ne\mu$, such that
$\gcd(A'_\lambda)=\gcd(A'_\mu)$. Denote this common gcd by $d$.
The displayed formula for $A'_\lambda-A'_\mu$ shows that
$d$ divides $a_{1n}$ (and $a_{n1}$).
It also divides all entries of $A'$.
Since $a_{1j}=0$ for $j<n$, it follows that $d$ divides
all entries of $A$. As $\gcd(A)=1$, we conclude that $d=1$.

We have shown that $\gcd(A'_\lambda)=1$ for some $\lambda\in F$.
By the induction hypothesis there exists $w\in R^{n-1}$
such that $f_{A'_\lambda}(w,w)=1$.
As $A$ and $A_\lambda$ are congruent,
there exists $v\in R^n$ such that $f_A(v,v)=1$.

\end{proof}

The second proposition is a skew-hermitian analog of the
first one.
We shall need the following definition.
Let $\nu_A$ denote the minimum degree of nonzero polynomials
$f_A(v,w)$ over all $v,w\in R^n$ with $f_A(v,v)=0$.

\begin{proposition} \label{sk-n}
If $A=(a_{ij})\in M_n(R)$, $A^*=-A$ and $\gcd(A)=1$,
then $A$ is congruent to the direct sum $B\oplus D$, where
\begin{equation} \label{mat-B}
B=\left(\begin{array}{cc}0&f\\-f^*&0\end{array}\right),
\end{equation}
with $f$ pure of degree $\nu_A$. Furthermore $ff^*$ divides 
all entries of $D$, i.e., $\det(B)$ is the second invariant
factor of $A$.
\end{proposition}

\begin{proof}
For the case $n=2$ see Proposition \ref{sk-2}. Thus let $n>2$.

After a suitable change of basis, we may assume that
$f_A(e_1,e_1)=0$ and that
there exists $w\in R^n$ such that $f_A(e_1,w)$ is nonzero and has
degree $\nu_A$. Thus $a_{11}=0$. By performing some
additional elementary congruence transformations, we may
also assume that $a_{12}\ne0$ has degree $\nu_A$
and that $a_{1j}=0$ for $j>2$.

Denote by $\sA$ the set of all skew-hermitian matrices
$X=(x_{ij})\in M_n(R)$ which are congruent to $A$ and
such that $x_{1j}=0$ for $j\ne2$ while $x_{12}$ has degree $\nu_A$.
For $X\in\sA$ let $d_X=\gcd(x_{12},x_{21},x_{22})$ where
we require $d_X$ to be monic.
Let $\sA_0$ denote the set of all $X\in\sA$ such that
$d_X$ has the minimum degree.
Without any loss of generality, we assume that
$A\in\sA_0$.

Our first objective is to show that $d_A$ is 1 or $t$.
Let $2\le r<s\le n$ and for $x\in R$ define $A_x\in\sA$ by
$$
A_x=(I_n+x^*E_{rs})A(I_n+xE_{sr})
$$
and set $d_x=d_{A_x}$. For $\lambda\in F$, the $(r,r)$-th entry
of $A_{\lambda x}$ is
\begin{equation} \label{rr}
a_{rr}+\lambda(a_{rs}x^*-a_{rs}^*x)+\lambda^2 a_{ss}xx^*.
\end{equation}

We take first $r=2$.
As $a_{12}$ has only finitely many monic divisors,
we can choose distinct $\alpha,\beta,\gamma\in F$ such that
$d_{\alpha x}=d_{\beta x}=d_{\gamma x}$.
Denote this common gcd by $d$. As the Vandermonde determinant
of $\alpha,\beta,\gamma$ is not 0, $d$ must
divide $a_{22}$, $a_{2s}x^*-a_{2s}^*x$ and $a_{ss}xx^*$.
It follows that $d$ divides $d_A$, and consequently
we must have $d=d_A$. By taking $x=1$, we infer
that $d_A$ divides the diagonal entries of $A$.
As $d_A$ divides $a_{2s}x^*-a_{2s}^*x$ for all $x\in R$,
we deduce that $d_A$ divides $ta_{2s}$.

Next we take $r>2$. Since $d_A$ must divide (\ref{rr})
for all $\lambda\in F$ and $x\in R$, we infer that 
$d_A$ divides $a_{rs}x^*-a_{rs}^*x$ for all $x\in R$.
Consequently, $d_A$ divides $ta_{rs}$. As $\gcd(A)=1$,
it follows that $d_A$ is either 1 or $t$.

We shall now rule out the possibility $d_A=t$.
Suppose that $d_A=t$. Assume that all entries $a_{2j}$
are divisible by $t$. In the above construction,
we take once again $r=2$ and choose $s>2$
such that $a_{sk}$ is not divisible by $t$
for some $k>1$ and $k\ne s$.
As above, we can choose a nonzero $\lambda\in F$ such that
$d_\lambda=d_A$. Then the $(2,k)$-th entry of $A_\lambda$
is not divisible by $t$.
Hence we can assume that one of the entries in the second
row, say $a_{23}$, is not divisible by $t$.

The $3\times3$ submatrix in the upper
left hand corner of $A$ has the form
\begin{equation*}
\left(\begin{array}{ccc}
0&at^k&0\\-a^*(-t)^k&bt&c\\0&-c^*&d\end{array}\right),
\end{equation*}
where $a,b,c$ are not divisible by $t$, $k\ge1$, $\gcd(a,a^*,b)=1$,
and the degree of $at^k$ is equal to $\nu_A$.
Moreover, by using Proposition \ref{her-2}, we may also
assume that $a$ is pure. Hence we can choose $x\in R$
such that $ax^*+a^*x=(-1)^kb$. Then the vector
$v=xe_1+t^{k-1}e_2$ is isotropic and
$$
f_A(v,e_1)=a^*t^{2k-1},\quad f_A(v,e_3)=c(-t)^{k-1}.
$$
Hence there exists $w\in R^n$ such that 
$$
f_A(v,w)=t^{k-1}\gcd(a^*,c),
$$
contradicting the fact that $at^k$ has degree $\nu_A$.
We conclude that $d_A=1$.

It is now easy to finish the proof.
By Proposition \ref{sk-2}, we may assume that the
$2\times2$ submatrix $B$ in the upper left hand corner
of $A$ has the form (\ref{mat-B}) with $f$ pure.
From the definition of $\nu_A$ it follows that the entries
$a_{1j}$, $j>2$, are divisible by $f$. By performing suitable
elementary congruence transformations, we may assume that 
all these entries are 0. A similar argument can be used
to make the entries $a_{2j}=0$ for $j>2$.

Thus we have $A=B\oplus D$ where $D=(d_{ij})$. Replace the zero
in the $(2,2)$ position of $A$ by $-d_{ii}$.
It follows from Proposition \ref{sk-2}
that this change can be achieved by a congruence
transformation on the block $B$.
Now add the $(i+2)$-nd row of $A$ to the second row and then
the $(i+2)$-nd column to the second column.
The entry in the $(2,2)$ position will become 0 again.
From the definition of $\nu_A$ it follows that the
$(2,j+2)$-nd entry of this new matrix must be divisible by $f^*$.
As this entry is equal to $d_{ij}$, we conclude
that all entries of $D$ are divisible by $f^*$. 
As $D^*=-D$, they are also divisible by $f$.
As $f$ is pure, all entries of $D$ are divisible by $ff^*$.

\end{proof}

We are now able to answer Kac's question.

\begin{theorem} \label{glavni}
If $A\in M_n(R)$ is hermitian or skew-hermitian,
then $A$ is congruent to the direct sum of $1\times1$
matrices and $2\times2$ matrices with zero diagonal.
\end{theorem}

\begin{proof}
As observed in section 2, we may assume that $\det(A)\ne0$
and $\gcd(A)=1$. We already know that the theorem is true
if $n\le2$. It remains to use induction and apply
the Propositions \ref{her-n} and \ref{sk-n}.
\end{proof}

To prove our second main result, we need the following
simple lemma.

\begin{lemma} \label{menjaj}
Let $A,B\in M_2(R)$ be skew-hermitian, $\gcd(A)=\gcd(B)=1$,
and $\det(A)=\det(B)\ne0$. Then $A$ and $B$ are congruent.
\end{lemma}

\begin{proof}
By Proposition \ref{sk-2}, we may assume that
$$
A=\left(\begin{array}{cc}0&ab\\-a^*b^*&0\end{array}\right), \quad
B=\left(\begin{array}{cc}0&ab^*\\-a^*b&0\end{array}\right),
$$
with $ab$ and $ab^*$ pure. There exist $x,y\in R_0$
such that $bb^*x-aa^*y=1$. If
$$
S=\left(\begin{array}{cc}b^*x&ay\\a^*&b\end{array}\right),
$$
then $S\in\SL_2(R)$ and $S^*BS=A$.
\end{proof}

Recall that two matrices $A,A'\in M_n(R)$ are said to be
{\em equivalent} if there exist $S,T\in\GL_n(R)$ such that
$A'=SAT$. A necessary and sufficient condition for $A$ and $A'$
to be equivalent is that they have the same invariant factors.

\begin{theorem} \label{invar}
Let $A,A'\in M_n(R)$ be both hermitian or both skew-hermitian.
If $A$ and $A'$ are equivalent, then they are congruent.
\end{theorem}

\begin{proof}
We use induction on $n$. The case $n=1$ is trivial. Let
$n>1$. Denote the invariant factors of $A$ (and $A'$) by
$f_1,\ldots,f_n$. If $f_n=0$ then we can use the induction
hypothesis. Assume that $f_n\ne0$. By dividing $A$ and
$A'$ by $f_1$, we may assume that $f_1=1$.

Now if $A$ and $A'$ are hermitian (resp. skew-hermitian)
then Proposition \ref{her-n} (resp. Proposition \ref{sk-n}
and Lemma \ref{menjaj})
allows us to finish the proof by using the induction hypothesis.

We shall give more details in the skew-hermitian case.
By Proposition \ref{sk-n} we may assume that $A=B\oplus D$,
where $B$ and $D$ are as stated there. Similarly, we may
assume that $A'=B'\oplus D'$. Since $\det(B)=f_2=\det(B')$,
Lemma \ref{menjaj} implies that $B$ and $B'$ are congruent.
Since $D$ and $D'$ have the same invariant factors,
they are congruent by the induction hypothesis. Hence
$A$ and $A'$ are congruent.

\end{proof}

\section{Canonical form under congruence}

In the next theorem we characterize the invariant factors of
hermitian and skew-hermitian matrices. Clearly these factors
have to be homogeneous.	

\begin{theorem}
Let $0\le r\le n$. Let
$f_1,\ldots,f_n$ be  a sequence of homogeneous elements in $R$
such that $f_1,\ldots,f_r$ are monic, each dividing the next one,
and $f_{r+1},\ldots,f_n$ are zero.
Then this sequence is the list of invariant factors
of a hermitian (resp. skew-hermitian) matrix $A\in M_n(R)$ of rank $r$
if and only if the following two conditions hold:
\begin{itemize}
\item[(i)] Any maximal subsequence $f_i,f_{i+1},\ldots,f_j$
consisting of consecutive nonzero odd (resp. even)
elements has even length. We shall write such subsequence as 
$$
g_i,h_i,g_{i+2},h_{i+2},\ldots,g_{j-1},h_{j-1}.
$$
\item[(ii)] For each $(g_k,h_k)$ as above,
$h_k=g_kp_kp_k^*$ with $p_k$ pure. 
\end{itemize}
\end{theorem}

\begin{proof}
We prove necessity by induction on $n$.
The cases $r=0$ and $n=1$ are trivial.
Let $r\ge1$ and $n\ge2$.
By replacing $A$ with $f_1^{-1}A$, we may assume that $f_1=1$.

If $A$ is hermitian, then Proposition \ref{her-n} shows that
$A$ is congruent to $(1)\oplus B$ and we can apply the induction
hypothesis to $B$ to finish the proof.

If $A$ is skew-hermitian, then $A$ is congruent
to the matrix $B\oplus D$ as stated in Proposition \ref{sk-n}.
In particular $f_2=\det(B)$ is even and not divisible by $t$.
We can now finish the proof by applying the induction hypothesis to $D$.

Sufficiency can be read off from the next theorem.
\end{proof}

It is now easy to obtain the canonical forms for hermitian
and skew-hermitian matrices under congruence.

\begin{theorem}
Let $A\in M_n(R)$ and $A^*=\ve A$, where $\ve=\pm$, 
let $r$ be the rank of $A$, and let $f_1,\ldots,f_n$ be the
invariant factors of $A$.
Form the direct sum, $B$, of the following blocks:
\begin{itemize}
\item[(i)] The $1\times1$ matrices $(f_i)$ for each $f_i$
such that $f_i^*=\ve f_i$.
\item[(ii)] The $2\times2$ matrices
$$
g_k\left(\begin{array}{cc}0&p_k\\\ve p_k^*&0\end{array}\right),
$$
for each pair $(g_k,h_k=g_kp_kp_k^*)$ constructed in the 
previous theorem.
\end{itemize}
Then $A$ is congruent to $B$. Moreover such $B$ is unique
up to the ordering of the diagonal blocks and the factorizations
$h_k=f_kp_kp_k^*$.
\end{theorem}

\begin{proof}
The matrices $A$ and $B$ have the same invariant factors.

\end{proof}

\section{Comments on other fields}

We introduce four conditions on a field $F$ assuming only that
the characteristic is not 2.

\begin{itemize}
\item[(K)] Kac's question has affirmative answer for the field $F$.
\item[(N)] The norm map $R\to R_0$ sending $x\to xx^*$ is onto.
\item[(U)] The quadratic form $x^2-ty^2$ over $R$ is universal.
\item[(I)] No element of $R_0$ is irreducible in $R$.
\end{itemize}

\begin{proposition}
For a field $F$ of characteristic $\ne2$, the above four
conditions are equivalent to each other.
\end{proposition}

\begin{proof}
$(K)\Rightarrow (N)$.
Let $\alpha\in F$, $\alpha\ne0$. As
$$
A=\left(\begin{array}{cc}-\alpha t&\alpha\\-\alpha&t\end{array}\right)
$$
is skew-hermitian but not diagonalizable, it must be congruent to
$$
\left(\begin{array}{cc}0&x\\-x^*&0\end{array}\right)
$$
for some $x\in R$.
Hence $\det(A)=-\alpha(t^2-\alpha)$ splits over $F$. We deduce
that $F$ is quadratically closed, i.e., it has no quadratic extensions.

It remains to show that if $a=1+t^2b$, with $b\in R_0$,
then $a=xx^*$ for some $x\in R$. This follows by applying
the above argument to
$$
\left(\begin{array}{cc}tb&1\\-1&t\end{array}\right).
$$

$(N)\Rightarrow (K)$. Our proofs are valid under this weaker hypothesis.

$(N)\Rightarrow (U)$. Let $\sigma:R\to R_0$ be the isomorphism
of $F$-algebras
sending $t$ to $t^2$. For $b\in R$ we have $\sigma(b)=zz^*$
for some $z\in R$. By writing $z=\sigma(x)+t\sigma(y)$,
$(x,y\in R)$, we obtain $\sigma(b)=\sigma(x)^2-t^2\sigma(y)^2$,
i.e., $b=x^2-ty^2$.

$(U)\Rightarrow (N)$. For $a\in R_0$ we have $a=\sigma(b)$
with $b\in R$. As $b=x^2-ty^2$ for some $x,y\in R$, we have
$a=zz^*$ with $z=\sigma(x)+t\sigma(y)$.

The equivalence between $(N)$ and $(I)$ is obvious.
\end{proof}

One can construct examples of fields $F$ satisfying the above conditions
without being algebraically closed. Start with a finite Galois
extension $K/E$ whose Galois group is not a 2-group.
Let $\sigma$ be an $E$-automorphism of an algebraic closure
$\overline{K}$ of $K$ whose restriction to $K$
is nontrivial and has odd order.
Then one can take $F$ to be the quadratic closure
of $\left(\overline{K}\right)^\sigma$.
In particular the quadratic closure of the prime field $F_p$
($p$ odd) is an example. On the other hand, it is easy to see
that the quadratic closure of the rationals does not satisfy
the condition (U).

In general, a hermitian or skew-hermitian matrix
$A\in M_n(R)$ need not be congruent to the direct sum of
{\em any} $1\times1$ or $2\times2$ matrices. For instance,
this is the case when $F$ is the real field and
$$
A=\left(\begin{array}{ccc}
t^2&1&0\\1&t^2&t\\0&-t&t^2\end{array}\right).
$$

\end{document}